\documentclass{article}
\usepackage{amssymb}
\usepackage{amsmath}
\usepackage{harvard}
\usepackage{latexsym}
\usepackage{amsmath,amsfonts,amssymb}
\usepackage{ifthen}
\usepackage{harvard}
\usepackage{latexsym}
\usepackage[a4paper,nohead]{geometry}

\provideboolean{uglystyle}
\setboolean{uglystyle}{true}
\newtheorem{theorem}{Theorem}[section]

\newcounter{tcnt}[theorem]

\newcounter{pcnt}[theorem]

\newcounter{ccnt}[theorem]

\newcounter{rcnt}[theorem]

\ifthenelse{\boolean{uglystyle}}{
\sloppy
\usepackage{leqno}
}{
\setlength{\parskip}{\medskipamount}
}
\input{tcilatex}
\fontsize{13pt}{13pt} \selectfont

\input epsf
\begin{document}

\title{Sparse Estimators and the Oracle Property, or the Return of Hodges'
Estimator}
\author{Hannes Leeb \\
Department of Statistics, Yale University\\
and \and Benedikt M. P\"{o}tscher \\
Department of Statistics, University of Vienna}
\date{First version: November 2004\\
This version: March 2007\\
}
\maketitle

\begin{abstract}
We point out some pitfalls related to the concept of an oracle property as
used in Fan and Li (2001, 2002, 2004) which are reminiscent of the
well-known pitfalls related to Hodges' estimator. The oracle property is
often a consequence of sparsity of an estimator. We show that any estimator
satisfying a sparsity property has maximal risk that converges to the
supremum of the loss function; in particular, the maximal risk diverges to
infinity whenever the loss function is unbounded. For ease of presentation
the result is set in the framework of a linear regression model, but
generalizes far beyond that setting. In a Monte Carlo study we also assess
the extent of the problem in finite samples for the smoothly clipped
absolute deviation (SCAD) estimator introduced in Fan and Li (2001). We find
that this estimator can perform rather poorly in finite samples and that its
worst-case performance relative to maximum likelihood deteriorates with
increasing sample size when the estimator is tuned to sparsity.

\textit{AMS 2000 Subject Classifications}: Primary 62J07, 62C99; secondary
62E20, 62F10, 62F12

\textit{Key words and phrases}: oracle property, sparsity, penalized maximum
likelihood, penalized least squares, Hodges' estimator, SCAD, Lasso, Bridge
estimator, hard-thresholding, maximal risk, maximal absolute bias,
non-uniform limits
\end{abstract}

\section{Introduction}

Recent years have seen an increased interest in penalized least squares and
penalized maximum likelihood estimation. Examples are the class of Bridge
estimators introduced by Frank and Friedman (1993), which includes
Lasso-type estimators as a special case (Knight and Fu (2000)), or the
smoothly clipped absolute deviation (SCAD) estimator introduced in Fan and
Li (2001) and further discussed in Fan and Li (2002, 2004), Fan and Peng
(2004), and Cai, Fan, Li, and Zhou (2005). As shown in Fan and Li (2001),
the SCAD estimator, with appropriate choice of the regularization (tuning)
parameter, possesses a sparsity property, i.e., it estimates zero components
of the true parameter vector exactly as zero with probability approaching
one as sample size increases while still being consistent for the non-zero
components. An immediate consequence of this sparsity property of the SCAD
estimator is that the asymptotic distribution of this estimator remains the
same whether or not the correct zero restrictions are imposed in the course
of the SCAD estimation procedure. [This simple phenomenon is true more
generally as pointed out, e.g., in P\"{o}tscher (1991, Lemma 1).] In other
words, with appropriate choice of the regularization parameter, the
asymptotic distribution of the SCAD estimator based on the overall model and
that of the SCAD estimator derived from the most parsimonious correct model
coincide. Fan and Li (2001) have dubbed this property the \textquotedblleft
oracle property\textquotedblright\ and have advertised this property of
their estimator.\footnote{%
The oracle property in the sense of Fan and Li should not be confused with
the notion of an oracle inequality as frequently used elsewhere in the
literature.} For appropriate choices of the regularization parameter, the
sparsity and the oracle property are also possessed by several -- but not
all -- members of the class of Bridge estimators (Knight and Fu (2000),
p.~1361, Zou (2006)). Similarly, suitably tuned thresholding procedures give
rise to sparse estimators.\footnote{%
These estimators do not satisfy the oracle property in case of
non-orthogonal design.} Finally, we note that traditional
post-model-selection estimators (e.g., maximum likelihood estimators
following model selection) based on a consistent model selection procedure
(for example, BIC or test procedures with suitably chosen critical values)
are another class of estimators that exhibit the sparsity and oracle
property; see P\"{o}tscher (1991) and Leeb and P\"{o}tscher (2005) for
further discussion. In a recent paper, Bunea (2004) uses such procedures in
a semiparametric framework and emphasizes the oracle property of the
resulting estimator; see also Bunea and McKeague (2005).

At first sight, the oracle property appears to be a desirable property of an
estimator as it seems to guarantee that, without knowing which components of
the true parameter are zero, we can do (asymptotically) as well as if we
knew the correct zero restrictions; that is, we can \textquotedblleft
adapt\textquotedblright\ to the unknown zero restrictions without paying a
price. This is too good to be true, and it is reminiscent of the
\textquotedblleft superefficiency\textquotedblright\ property of the Hodges'
estimator; and justly so, since Hodges' estimator in its simplest form is a
hard-thresholding estimator exhibiting the sparsity and oracle property.
[Recall that in its simplest form Hodges' estimator for the mean of an $%
N(\mu ,1)$-distribution is given by the arithmetic mean $\bar{y}$ of the
random sample of size $n$ if $\left\vert \bar{y}\right\vert $ exceeds the
threshold $n^{-1/4}$, and is given by zero otherwise.] Now, as is
well-known, e.g., from Hodges' example, the oracle property is an asymptotic
feature that holds only \textit{pointwise} in the parameter space and gives
a misleading picture of the actual finite-sample performance of the
estimator. In fact, the finite sample properties of an estimator enjoying
the oracle property are often markedly different from what the pointwise
asymptotic theory predicts; e.g., the finite sample distribution can be
bimodal regardless of sample size, although the pointwise asymptotic
distribution is normal. This is again well-known for Hodges' estimator. For
a more general class of post-model-selection estimators possessing the
sparsity and the oracle property this is discussed in detail in Leeb and P%
\"{o}tscher (2005), where it is, e.g., also shown that the finite sample
distribution can \textquotedblleft escape to infinity\textquotedblright\
along appropriate local alternatives although the pointwise asymptotic
distribution is perfectly normal.\footnote{%
That pointwise asymptotics can be misleading in the context of model
selection has been noted earlier in Hosoya (1984), Shibata (1986a), P\"{o}%
tscher (1991), and Kabaila (1995, 2002).} See also Knight and Fu (2000,
Section 3) for related results for Bridge estimators. Furthermore,
estimators possessing the oracle property are certainly not exempt from the
Hajek-LeCam local asymptotic minimax theorem, further eroding support for
the belief that these estimators are as good as the \textquotedblleft
oracle\textquotedblright\ itself (i.e., the infeasible \textquotedblleft
estimator\textquotedblright\ that uses the information which components of
the parameter are zero).

The above discussion shows that the reasoning underlying the oracle property
is misguided. Even worse, estimators possessing the sparsity property (which
often entails the oracle property) necessarily have dismal finite sample
performance: It is well-known for Hodges' estimator that the maximal
(scaled) mean squared error grows without bound as sample size increases
(e.g., Lehmann and Casella (1998), p.442), whereas the standard maximum
likelihood estimator has constant finite quadratic risk. In this note we
show that a similar unbounded risk result is in fact true for \emph{any}
estimator possessing the sparsity property. This means that there is a
substantial price to be paid for sparsity even though the oracle property
(misleadingly) seems to suggest otherwise. As discussed in more detail
below, the bad risk behavior is a \textquotedblleft local\textquotedblright\
phenomenon and furthermore occurs at points in the parameter space that are
\textquotedblleft sparse\textquotedblright\ in the sense that some of their
coordinates are equal to zero. For simplicity of presentation and for
reasons of comparability with the literature cited earlier, the result will
be set in the framework of a linear regression model, but inspection of the
proof shows that it easily extends far beyond that framework. For related
results in the context of traditional post-model-selection estimators see
Yang (2005) and Leeb and P\"{o}tscher (2005, Appendix C);\footnote{%
The unboundedness of the maximal (scaled) mean squared error of estimators
following BIC-type model selection has also been noted in Hosoya (1984),
Shibata (1986b), and Foster and George (1994).} cf.~also the discussion on
\textquotedblleft partially\textquotedblright\ sparse estimators towards the
end of Section 2 below. The theoretical results in Section 2 are rounded out
by a Monte Carlo study in Section 3 that demonstrates the extent of the
problem in finite samples for the SCAD estimator of Fan and Li (2001). The
reasons for concentrating on the SCAD estimator in the Monte Carlo study are
(i) that the finite-sample risk behavior of traditional post-model-selection
estimators is well-understood (Judge and Bock (1978), Leeb and P\"{o}tscher
(2005)) and (ii) that the SCAD estimator -- especially when tuned to
sparsity -- has been highly advertised as a superior procedure in Fan and Li
(2001) and subsequent papers mentioned above.

\section{Bad Risk Behavior of Sparse Estimators}

Consider the linear regression model 
\begin{equation}
y_{t}\quad =\quad x_{t}^{\prime }\theta +\epsilon _{t}\qquad (1\leq t\leq n)
\label{model}
\end{equation}%
where the $k\times 1$ nonstochastic regressors $x_{t}$ satisfy $%
n^{-1}\sum_{t=1}^{n}x_{t}x_{t}^{\prime }$ $\rightarrow $ $Q>0$ as $%
n\rightarrow \infty $ and the prime denotes transposition. The errors $%
\epsilon _{t}$ are assumed to be independent identically distributed with
mean zero and finite variance $\sigma ^{2}$. Without loss of generality we
freeze the variance at $\sigma ^{2}=1$.\footnote{%
If the variance is not frozen at $\sigma ^{2}=1$, the results below
obviously continue to hold for each fixed value of $\sigma ^{2}$, and hence
hold a fortiori if the supremum in (\ref{2})--(\ref{3}) below\ is also taken
over $\sigma ^{2}$.} Furthermore, we assume that $\epsilon _{t}$ has a
density $f$ that possesses an absolutely continuous derivative $df/dx$
satisfying%
\begin{equation*}
0<\dint\limits_{-\infty }^{\infty }\left( (df(x)/dx)/f(x)\right)
^{2}f(x)dx<\infty .
\end{equation*}%
Note that the conditions on $f$ guarantee that the information of $f$ is
finite and positive. These conditions are obviously satisfied in the special
case of normally distributed errors. Let $P_{n,\theta }$ denote the
distribution of the sample $(y_{1},\ldots ,y_{n})^{\prime }$ and let $%
E_{n,\theta }$ denote the corresponding expectation operator. For $\theta
\in \mathbb{R}^{k}$, let $r(\theta )$ denote a $k\times 1$ vector with
components $r_{i}(\theta )$ where $r_{i}(\theta )=0$ if $\theta _{i}=0$ and $%
r_{i}(\theta )=1$ if $\theta _{i}\neq 0$. An estimator $\hat{\theta}$ for $%
\theta $ based on the sample $(y_{1},\ldots ,y_{n})^{\prime }$ is said to
satisfy the sparsity-type condition if for every $\theta \in \mathbb{R}^{k}$%
\begin{equation}
P_{n,\theta }\left( r(\hat{\theta})\leq r(\theta )\right) \rightarrow 1
\label{1}
\end{equation}%
holds for $n\rightarrow \infty $, where the inequality sign is to be
interpreted componentwise. That is, the estimator is guaranteed to find the
zero components of $\theta $ with probability approaching one as $%
n\rightarrow \infty $. Clearly, any sparse estimator satisfies (\ref{1}). In
particular, the SCAD estimator as well as certain members of the class of
Bridge estimators satisfy (\ref{1}) for suitable choices of the
regularization parameter as mentioned earlier. Also, any
post-model-selection estimator based on a consistent model selection
procedure clearly satisfies (\ref{1}). All these estimators are additionally
also consistent for $\theta $, and hence in fact satisfy the stronger
condition $P_{n,\theta }(r(\hat{\theta})=r(\theta ))\rightarrow 1$ for all $%
\theta \in \mathbb{R}^{k}$. [Condition (\ref{1}) by itself is of course also
satisfied by nonsensical estimators like $\hat{\theta}\equiv 0$, but is all
that is needed to establish the subsequent result.] We now show that any
estimator satisfying the sparsity-type condition (\ref{1}) has quite bad
finite sample risk properties. For purposes of comparison we note that the
(scaled) mean squared error of the least squares estimator $\hat{\theta}%
_{LS} $\ satisfies%
\begin{equation*}
E_{n,\theta }\left[ n(\hat{\theta}_{LS}-\theta )^{\prime }(\hat{\theta}%
_{LS}-\theta )\right] =\limfunc{trace}\left[ \left(
n^{-1}\sum_{t=1}^{n}x_{t}x_{t}^{\prime }\right) ^{-1}\right]
\end{equation*}%
which converges to $\limfunc{trace}(Q^{-1})$, and thus remains bounded as
sample size increases.

\begin{theorem}
\footnote{%
Theorem 2.1 and the ensuing discussion continue to apply if the regressors $%
x_{t}$ as well as the errors $\epsilon _{t}$ are allowed to depend on sample
size $n$, at least if the errors are normally distributed. The proof is
analogous, except that one uses direct computation and LeCam's first lemma
(cf., e.g., Lemma A.1 in Leeb and P\"{o}tscher (2006)) instead of Koul and
Wang (1984) to verify contiguity. Also, the results continue to hold if the
design matrix satisfies $\delta _{n}^{-1}\sum_{t=1}^{n}x_{t}x_{t}^{\prime }$ 
$\rightarrow $ $Q>0$ for some positive sequence $\delta _{n}$ other than $n$%
, provided that the scaling factor $n^{1/2}$ is replaced by $\delta
_{n}^{1/2}$ throughout.}Let $\hat{\theta}$ be an arbitrary estimator for $%
\theta $ that satisfies the sparsity-type condition (\ref{1}). Then the
maximal (scaled) mean squared error of $\hat{\theta}$ diverges to infinity
as $n\rightarrow \infty $, i.e., 
\begin{equation}
\sup_{\theta \in \mathbb{R}^{k}}E_{n,\theta }\left[ n(\hat{\theta}-\theta
)^{\prime }(\hat{\theta}-\theta )\right] \rightarrow \infty   \label{2}
\end{equation}%
for $n\rightarrow \infty $. More generally, let $l:\mathbb{R}^{k}\rightarrow 
\mathbb{R}$ be a nonnegative loss function. Then 
\begin{equation}
\sup_{\theta \in \mathbb{R}^{k}}E_{n,\theta }l(n^{1/2}(\hat{\theta}-\theta
))\rightarrow \sup_{s\in \mathbb{R}^{k}}l(s)  \label{3}
\end{equation}%
for $n\rightarrow \infty $. In particular, if the loss function $l$ is
unbounded then the maximal risk associated with $l$ diverges to infinity as $%
n\rightarrow \infty $.
\end{theorem}

The theorem says that, whatever the loss function, the maximal risk of a
sparse estimator is -- in large samples -- as bad as it possibly can be.

Upon choosing $l(s)=\left\vert s_{i}\right\vert $, where $s_{i}$ denotes the 
$i$-th coordinate of $s$, relation (\ref{3}) shows that also the maximal
(scaled) absolute bias of each component $\hat{\theta}_{i}$ diverges to
infinity.

Applying relation (\ref{3}) to the loss function $l^{\ast }(s)=l(c^{\prime
}s)$ shows that (\ref{3}) holds mutatis mutandis also for estimators $%
c^{\prime }\hat{\theta}$ of arbitrary linear contrasts $c^{\prime }\theta $.
In particular, using quadratic loss $l^{\ast }(s)=(c^{\prime }s)^{2}$, it
follows that also the maximal (scaled) mean squared error of the linear
contrast $c^{\prime }\hat{\theta}$ goes to infinity as sample size
increases, provided $c\neq 0$.

\textbf{Proof of Theorem 2.1: }It suffices to prove (\ref{3}).\footnote{%
Note that the expectations in (\ref{2}) and (\ref{3}) are always
well-defined.} Now, with $\theta _{n}=-n^{-1/2}s$, $s\in \mathbb{R}^{k}$
arbitrary, we have%
\begin{align}
\sup_{u\in \mathbb{R}^{k}}l(u)& \geq \sup_{\theta \in \mathbb{R}%
^{k}}E_{n,\theta }l(n^{1/2}(\hat{\theta}-\theta ))\geq E_{n,\theta
_{n}}l(n^{1/2}(\hat{\theta}-\theta _{n}))  \notag \\
& \geq E_{n,\theta _{n}}[l(n^{1/2}(\hat{\theta}-\theta _{n}))\boldsymbol{1}(%
\hat{\theta}=0)]=l(-n^{1/2}\theta _{n})P_{n,\theta _{n}}(r(\hat{\theta})=0) 
\notag \\
& =l(s)P_{n,\theta _{n}}(r(\hat{\theta})=0).  \label{4}
\end{align}%
By the sparsity-type condition we have that $P_{n,0}(r(\hat{\theta}%
)=0)\rightarrow 1$ as $n\rightarrow \infty $. Since the model is locally
asymptotically normal under our assumptions (Koul and Wang (1984), Theorem
2.1 and Remark 1; Hajek and Sidak (1967), p.213), the sequence of
probability measures $P_{n,\theta _{n}}$ is contiguous w.r.t. the sequence $%
P_{n,0}$. Consequently, the far r.h.s. of (\ref{4}) converges to $l(s)$.
Since $s\in \mathbb{R}^{k}$ was arbitrary, the proof is complete.{\small %
\hfill }$\blacksquare $

Inspection of the proof shows that Theorem 2.1 remains true if the supremum
of the risk in (\ref{3}) is taken only over open balls of radius $\rho _{n}$
centered at the origin as long as $n^{1/2}\rho _{n}\rightarrow \infty $.
Hence, the bad risk behavior is a local phenomenon that occurs in a part of
the parameter space where one perhaps would have expected the largest gain
over the least squares estimator due to the sparsity property. [If the
supremum of the risk in (\ref{3}) is taken over the open balls of radius $%
n^{-1/2}\rho $ centered at the origin where $\rho >0$ is now fixed, then the
proof still shows that the limit inferior of this supremum is not less than $%
\sup_{\left\Vert s\right\Vert <\rho }l(s)$.] Furthermore, for quadratic loss 
$l(s)=s^{\prime }s$, a small variation of the proof shows that these
\textquotedblleft local\textquotedblright\ results continue to hold if the
open balls over which the supremum is taken are not centered at the origin,
but at an arbitrary $\theta $, as long as $\theta $ possesses at least one
zero component. [It is easy to see that this is more generally true for any
nonnegative loss function $l$ satisfying, e.g., $l(s)\geq l(\pi _{i}(s))$
for every $s\in \mathbb{R}^{k}$ and an index $i$ with $\theta _{i}=0$, where 
$\pi _{i}$ represents the projection on the $i$-th coordinate axis.]

Inspection of the proof also shows that -- at least in the case of quadratic
loss -- the element $s$ can be chosen to point in the direction of a
standard basis vector. This then shows that the bad risk behavior occurs at
parameter values that themselves are \textquotedblleft
sparse\textquotedblright\ in the sense of having many zero coordinates.

If the quadratic loss $n(\hat{\theta}-\theta )^{\prime }(\hat{\theta}-\theta
)$\ in (\ref{2}) is replaced by the weighted quadratic loss $(\hat{\theta}%
-\theta )^{\prime }\sum_{t=1}^{n}x_{t}x_{t}^{\prime }(\hat{\theta}-\theta )$%
, then the corresponding maximal risk again diverges to infinity. More
generally, let $l_{n}$ be a nonnegative loss function that may depend on
sample size. Inspection of the proof of Theorem 2.1 shows that%
\begin{equation}
\limsup_{n\rightarrow \infty }\sup_{u\in \mathbb{R}^{k}}l_{n}(u)\geq
\limsup_{n\rightarrow \infty }\sup_{\left\Vert \theta \right\Vert
<n^{-1/2}\rho }E_{n,\theta }l_{n}(n^{1/2}(\hat{\theta}-\theta ))\geq
\sup_{\left\Vert u\right\Vert <\rho }\limsup_{n\rightarrow \infty }l_{n}(u),
\label{5}
\end{equation}%
\begin{equation}
\liminf_{n\rightarrow \infty }\sup_{u\in \mathbb{R}^{k}}l_{n}(u)\geq
\liminf_{n\rightarrow \infty }\sup_{\left\Vert \theta \right\Vert
<n^{-1/2}\rho }E_{n,\theta }l_{n}(n^{1/2}(\hat{\theta}-\theta ))\geq
\sup_{\left\Vert u\right\Vert <\rho }\liminf_{n\rightarrow \infty }l_{n}(u)
\label{6}
\end{equation}%
hold for any $0<\rho \leq \infty $. [In case $0<\rho <\infty $, the lower
bounds in (\ref{5})-(\ref{6}) can even be improved to $\limsup_{n\rightarrow
\infty }\sup_{\left\Vert u\right\Vert <\rho }l_{n}(u)$ and $%
\liminf_{n\rightarrow \infty }\sup_{\left\Vert u\right\Vert <\rho }l_{n}(u)$%
, respectively.\footnote{%
Note that the local asymptotic normality condition in Koul and Wang (1984)
as well as the result in Lemma A.1 in Leeb and P\"{o}tscher (2006) imply
contiguity of $P_{n,\theta _{n}}$ and $P_{n,0}$ not only for $\theta
_{n}=\gamma /n^{1/2}$ but more generally for $\theta _{n}=\gamma
_{n}/n^{1/2} $ with $\gamma _{n}$ a bounded sequence.} It then follows that
in case $\rho =\infty $ the lower bounds in (\ref{5})-(\ref{6}) can be
improved to $\sup_{0<\tau <\infty }\limsup_{n\rightarrow \infty
}\sup_{\left\Vert u\right\Vert <\tau }l_{n}(u)$ and $\sup_{0<\tau <\infty
}\liminf_{n\rightarrow \infty }\sup_{\left\Vert u\right\Vert <\tau }l_{n}(u)$%
, respectively.]

Next we briefly discuss the case where an estimator $\hat{\theta}$ only has
a \textquotedblleft partial\textquotedblright\ sparsity property (and
consequently a commensurable oracle property) in the following sense:
Suppose the parameter vector $\theta $ is partitioned as $\theta =(\alpha
^{\prime },\beta ^{\prime })^{\prime }$ and the estimator $\hat{\theta}=(%
\hat{\alpha}^{\prime },\hat{\beta}^{\prime })^{\prime }$ only finds the true
zero components in the subvector $\beta $ with probability converging to
one. E.g., $\hat{\theta}$ is a traditional post-model-selection estimator
based on a consistent model selection procedure that is designed to only
identify the zero components in $\beta $. A minor variation of the proof of
Theorem 2.1 immediately shows again that the maximal (scaled) mean squared
error of $\hat{\beta}$, and hence also of $\hat{\theta}$, diverges to
infinity, and the same is true for linear combinations $d^{\prime }\hat{\beta%
}$ as long as $d\neq 0$. [This immediately extends to linear combinations $%
c^{\prime }\hat{\theta}$, as long as $c$ charges at least one coordinate of $%
\hat{\beta}$ with a nonzero coefficient.]\footnote{%
In fact, this variation of the proof of Theorem 2.1 shows that the supremum
of $E_{n,\theta }l(n^{1/2}(\hat{\beta}-\beta ))$, where $l$ is an arbitrary
nonegative loss function, again converges to the supremum of the loss
function.} However, if the parameter of interest is $\alpha $ rather than $%
\beta $, Theorem 2.1 and its proof (or simple variations thereof) do not
apply to the mean squared error of $\hat{\alpha}$ (or its linear contrasts).
Nevertheless, the maximal (scaled) mean squared error of $\hat{\alpha}$ can
again be shown to diverge to infinity, at least for traditional
post-model-selection estimators $\hat{\theta}$ based on a consistent model
selection procedure; see Leeb and P\"{o}tscher (2005, Appendix C).

While the above results are set in the framework of a linear regression
model with nonstochastic regressors, it is obvious from the proof that they
extend to much more general models such as regression models with stochastic
regressors, semiparametric models, nonlinear models, time series models,
etc., as long as the contiguity property used in the proof is satisfied.
This is in particular the case whenever the model is locally asymptotically
normal, which in turn is typically the case under standard regularity
conditions for maximum likelihood estimation.

\section{Numerical Results on the Finite Sample Performance of the SCAD
Estimator}

We replicate and extend Monte Carlo simulations of the performance of the
SCAD estimator given in Example 4.1 of Fan and Li (2001); we demonstrate
that this estimator, when tuned to enjoy a sparsity property and an oracle
property, can perform quite unfavorably in finite samples. Even when not
tuned to sparsity, we show that the SCAD estimator can perform worse than
the least squares estimator in parts of the parameter space, something that
is not brought out in the simulation study in Fan and Li (2001) as they
conducted their simulation only at a single point in the parameter space
(which happens to be favorable to their estimator).

Consider $n$ independent observations from the linear model (1) with $k=8$
regressors, where the errors $\epsilon _{t}$ are standard normal and are
distributed independently of the regressors. The regressors $x_{t}$ are
assumed to be multivariate normal with mean zero. The variance of each
component of $x_{t}$ is equal to $1$ and the correlation between the $i$-th
and the $j$-th component of $x_{t}$, i.e., $x_{t,i}$ and $x_{t,j}$, is $\rho
^{|i-j|}$ with $\rho =0.5$. Fan and Li (2001) consider this model with $n=40$%
, $n=60$, and with the true parameter equal to $\theta
_{0}=(3,1.5,0,0,2,0,0,0)^{\prime }$; cf. also Tibshirani (1996, Section
7.2). We consider a whole range of true values for $\theta $ at various
sample sizes, namely $\theta _{n}=\theta _{0}+(\gamma /\sqrt{n})\times \eta $
for some vector $\eta $ and for a range of $\gamma $'s as described below.
We do this because (i) considering only one choice for the true parameter in
a simulation may give a wrong impression of the actual performance of the
estimators considered, and (ii) because our results in Section 2 suggest
that the risk of sparse estimators can be large for parameter vectors which
have some of its components close to, but different from, zero.

The SCAD estimator is defined as a solution to the problem of minimizing the
penalized least squares objective function%
\begin{equation*}
\frac{1}{2}\sum_{t=1}^{n}(y_{t}-x_{t}^{\prime }\theta
)^{2}+n\sum_{i=1}^{k}p_{\lambda }(\left\vert \theta _{i}\right\vert )
\end{equation*}%
where the penalty function $p_{\lambda }$ is defined in Fan and Li (2001)
and $\lambda \geq 0$ is a tuning parameter. The penalty function $p_{\lambda
}$ contains also another tuning parameter $a$, which is set equal to 3.7
here, resulting in a particular instance of the SCAD estimator which is
denoted by SCAD2 in Example 4.1 of Fan and Li (2001). According to Theorem 2
in Fan and Li (2001) the SCAD estimator is guaranteed to satisfy the
sparsity property if $\lambda \rightarrow 0$ and $\sqrt{n}\lambda
\rightarrow \infty $ as samples size $n$ goes to infinity.

Using the MATLAB code provided to us by Runze Li, we have implemented the
SCAD2 estimator in R. [The code is available from the first author on
request.] Two types of performance measures are considered: The `median
relative model error' studied by Fan and Li (2001), and the relative mean
squared error. The median relative model error is defined as follows: For an
estimator $\hat{\theta}$ for $\theta $, define the model error $ME(\hat{%
\theta})$ by $ME(\hat{\theta})=(\hat{\theta}-\theta )^{\prime }\Sigma (\hat{%
\theta}-\theta )$, where $\Sigma $ denotes the variance/covariance matrix of
the regressors. Now define the relative model error of $\hat{\theta}$
(relative to least squares) by $ME(\hat{\theta})/ME(\hat{\theta}_{LS})$,
with $\hat{\theta}_{LS}$ denoting the least squares estimator based on the
overall model. The median relative model error is then given by the median
of the relative model error. The relative mean squared error of $\hat{\theta}
$ is given by $E[(\hat{\theta}-\theta )^{\prime }(\hat{\theta}-\theta )]/E[(%
\hat{\theta}_{LS}-\theta )^{\prime }(\hat{\theta}_{LS}-\theta )]$.\footnote{%
The mean squared error of \ $\hat{\theta}_{LS}$ is given by $E\limfunc{trace}%
((X^{\prime }X)^{-1})$ which equals $\limfunc{trace}(\Sigma
^{-1})/(n-9)=38/(3n-27)$ by von Rosen (1988, Theorem 3.1).} Note that we
have scaled the performance measures such that both of them are identical to
unity for $\hat{\theta}=\hat{\theta}_{LS}$.

\noindent \textbf{Setup I:} For SCAD2 the tuning parameter $\lambda $ is
chosen by generalized cross-validation (cf. Section 4.2 of Fan and Li
(2001)). In the original study in Fan and Li (2001), the range of $\lambda $%
's considered for generalized cross-validation at sample sizes $n=40$ and $%
n=60$ is $\{\delta (\hat{\sigma}/\sqrt{n}):\;\delta =0.9,1.1,1.3,\dots ,2\}$%
; here, $\hat{\sigma}^{2}$ denotes the usual unbiased variance estimator
obtained from a least-squares fit of the overall model. For the simulations
under Setup I, we re-scale this range of $\lambda $'s by $\log {n}/\log {60}$%
. With this, our results for $\gamma =0$ replicate those in Fan and Li
(2001) for $n=60$; for the other (larger) sample sizes that we consider, the
re-scaling guarantees that $\lambda \rightarrow 0$ and $\sqrt{n}\lambda
\rightarrow \infty $ and hence, in view of Theorem 2 in Fan and Li (2001),
guarantees that the resulting estimator enjoys the sparsity condition. [For
another choice of $\lambda $ see Setup VI.] We compute Monte Carlo estimates
for both the median relative model error and the relative mean squared error
of the SCAD2 estimator for a range of true parameter values, namely $\theta
_{n}=\theta _{0}+(\gamma /\sqrt{n})\times (0,0,1,1,0,1,1,1)^{\prime }$ for
101 equidistant values of $\gamma $ between $0$ and $8$, and for sample
sizes $n=60$, $120$, $240$, $480$, and $960$, each based on 500 Monte Carlo
replications (for comparison, Fan and Li (2001) use 100 replications). Note
that the performance measures are symmetric about $\gamma =0$, and hence are
only reported for nonnegative values of $\gamma $. The results are
summarized in Figure~1 below. [For better readability, points in Figure~1
are joined by lines.]

\begin{center}
\begin{tabular}{cc}
\epsfxsize=7cm\epsfbox{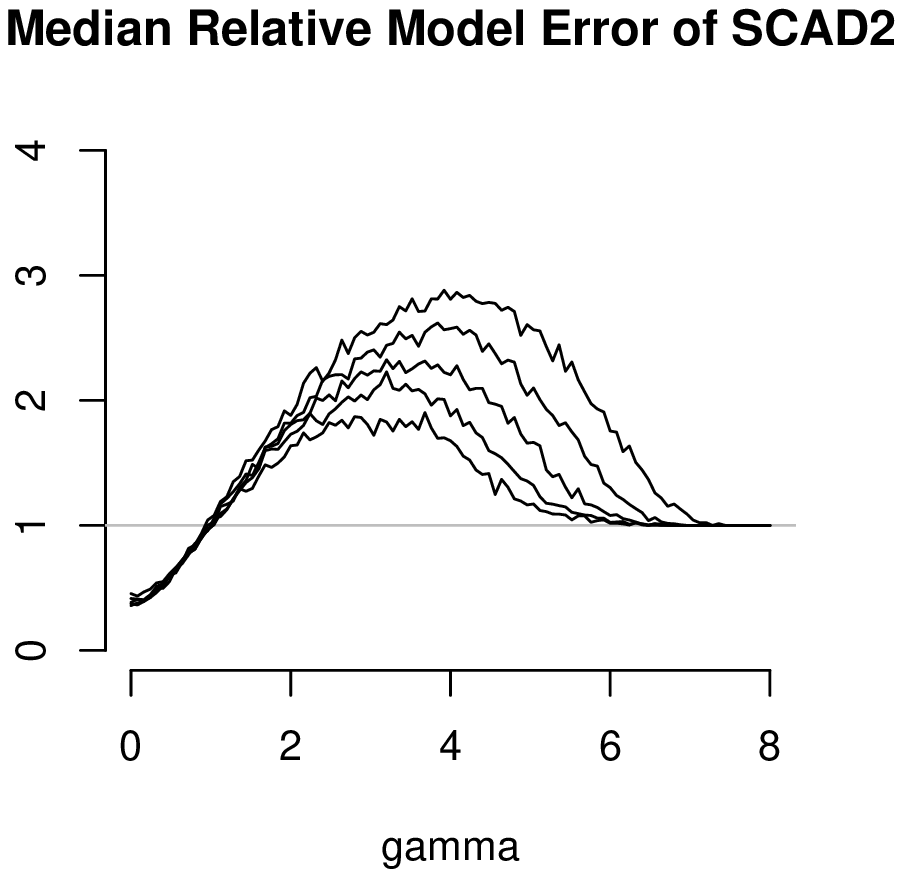} & \epsfxsize=7cm\epsfbox{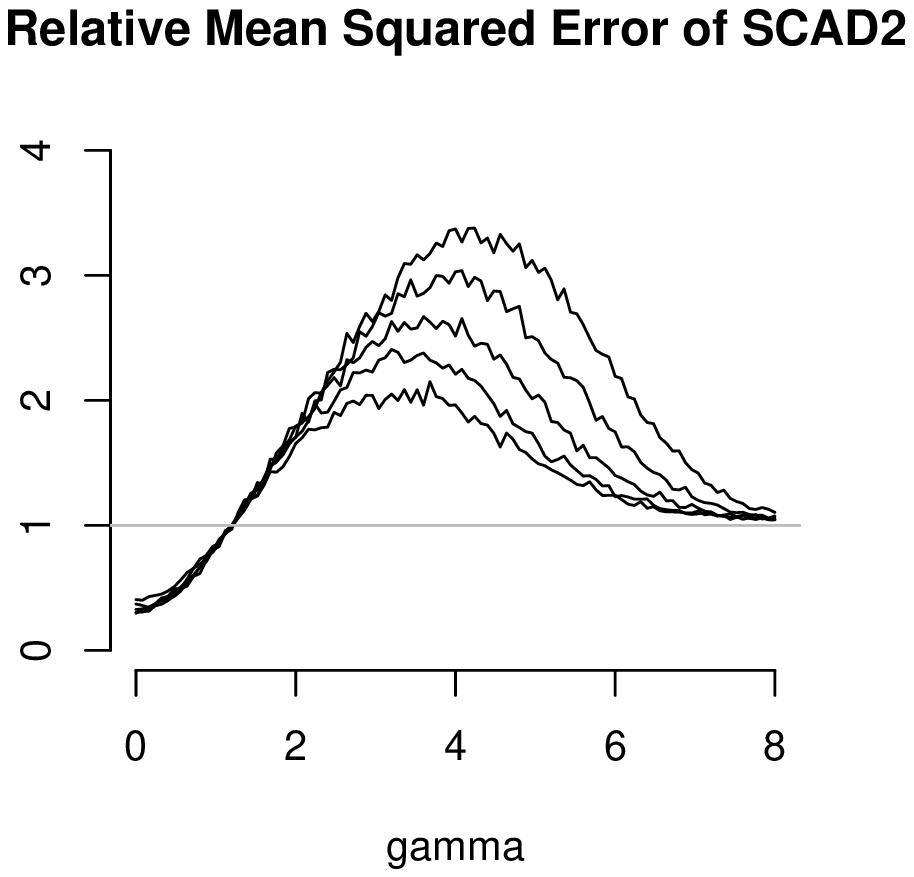}%
\end{tabular}
\end{center}

\begin{quote}
{\small Figure~1: Monte Carlo performance estimates under the true parameter 
$\theta _{n}=\theta _{0}+(\gamma /\sqrt{n})\times (0,0,1,1,0,1,1,1)^{\prime
} $ , as a function of $\gamma $. The left panel gives the estimated median
relative model error of SCAD2 for sample sizes $n=60,120,240,480,960$. The
right panel gives the corresponding results for the estimated relative mean
squared error of SCAD2. Larger sample sizes correspond to larger maximal
errors. For comparison, the gray line at one indicates the performance of
the ordinary least squares estimator.}
\end{quote}

In the Monte Carlo study of Fan and Li (2001), only the parameter value $%
\theta _{0}$ is considered. This corresponds to the point $\gamma =0$ in the
panels of Figure~1. At that particular point in the parameter space, SCAD2
compares quite favorably with least squares. However, Figure~1 shows that
there is a large range of parameters where the situation is reversed. In
particular, we see that SCAD2 can perform quite unfavorably when compared to
least squares if the true parameter, i.e., $\theta _{n}$, is such that some
of its components are close to, but different from, zero. In line with
Theorem~2.1, we also see that the worst-case performance of SCAD2
deteriorates with increasing sample size: For $n=60$, ordinary least squares
beats SCAD2 in terms of worst-case performance by a factor of about 2 in
both panels of Figure~1; for $n=960$, this factor has increased to about 3;
and increasing $n$ further makes this phenomenon even more pronounced. We
also see that, for increasing $n$, the location of the peak moves to the
right in Figure~1, suggesting that the worst-case performance of SCAD2
(among parameters of the form $\theta _{n}=(\gamma /\sqrt{n})\times
(0,0,1,1,0,1,1,1)^{\prime }$) is attained at a value $\gamma _{n}$, which is
such that $\gamma _{n}\rightarrow \infty $ with $n$. In view of the proof of
Theorem~2.1, this is no surprise.\footnote{%
See Section 2.1 and Footnote 14 in Leeb and P\"{o}tscher (2005) for related
discussion.} [Of course, there may be other parameters at any given sample
size for which SCAD2 performs even worse.] Our simulations thus demonstrate:
If each component of the true parameter is either very close to zero or
quite large (where the components' size has to be measured relative to
sample size), then the SCAD estimator performs well. However, if some
component is in-between these two extremes, then the SCAD estimator performs
poorly. In particular, the estimator can perform poorly precisely in the
important situation where it is statistically difficult to decide whether
some component of the true parameter is zero or not. Poor performance is
obtained in the worst case over a neighborhood of one of the
lower-dimensional models, where the `diameter' of the neighborhood goes to
zero slower than $1/\sqrt{n}$.

We have also re-run our simulations for other experimental setups; the
details are given below. Since our findings for these other setups are
essentially similar to those summarized in Figure~1, we first give a brief
overview of the other setups and summarize the results before proceeding to
the details. In Setups II and III we consider slices of the $8$-dimensional
performance measure surfaces corresponding to directions other than the one
used in Setup I: In Setup~II the true parameter is of the form $\theta
_{0}+(\gamma /\sqrt{n})\times (0,0,1,1,0,0,0,0)^{\prime }$, i.e., we
consider the case where some components are exactly zero, some are large,
and others are in-between. In Setup~III, we consider a scenario in-between
Setup~I and Setup~II, namely the case where the true parameter is of the
form $\theta _{0}+(\gamma /\sqrt{n})\times
(0,0,1,1,0,1/10,1/10,1/10)^{\prime }$. The method for choosing $\lambda $ in
these two setups is the same as in Setup I. The results in these additional
setups are qualitatively similar to those shown in Figure~1 but slightly
less pronounced. In further setups we also consider various other rates for
the SCAD tuning parameter $\lambda $. By Theorem 2 of Fan and Li (2001), the
SCAD estimator is sparse if $\lambda \rightarrow 0$ and $\sqrt{n}\lambda
\rightarrow \infty $; as noted before, for Figure~1, $\lambda $ is chosen by
generalized cross-validation from the set $\Lambda _{n}=\{\delta (\hat{\sigma%
}/\sqrt{n})(\log (n)/\log (60)):\;\delta =0.9,1.1,1.3,\dots ,2\}$; i.e., we
have $\sqrt{n}\lambda =O_{p}(\log (n))$. The magnitude of $\lambda $ has a
strong impact on the performance of the estimator. Smaller values result in
`less sparse' estimates, leading to less favorable performance relative to
least squares at $\gamma =0$, but at the same time leading to less
unfavorable worst-case performance; the resulting performance curves are
`flatter' than those in Figure~1. Larger values of $\lambda $ result in
`more sparse' estimates, improved performance at $\gamma =0$, and more
unfavorable worst-case performance; this leads to performance curves that
are `more spiked' than those in Figure~1. In Setups~\ IV and V we have
re-run our simulations with $\gamma $ chosen from a set $\Lambda _{n}$ as
above, but with $\log (n)/\log (60)$ replaced by $(n/60)^{1/10}$ as well as
by $(n/60)^{1/4}$, resulting in $\sqrt{n}\lambda =O_{p}(n^{1/10})$ and $%
\sqrt{n}\lambda =O_{p}(n^{1/4})$, respectively. In Setup~IV, where $\sqrt{n}%
\lambda =O_{p}(n^{1/10})$, we get results similar to, but less pronounced
than, Figure~1; this is because Setup~IV leads to $\lambda $'s smaller than
in Setup~I. In Setup~V, where $\sqrt{n}\lambda =O_{p}(n^{1/4})$, we get
similar but more pronounced results when compared to Figure 1; again, this
is so because Setup~V leads to larger $\lambda $'s than Setup~I. A final
setup (Setup VI) in which we do not enforce the conditions for sparsity is
discussed below after the details for Setups II-V are presented.

\noindent \textbf{Setups~II and~III:} In Setup~II, we perform the same Monte
Carlo study as in Setup~I, the only difference being that the range of $%
\theta $'s is now $\theta _{n}=\theta _{0}+(\gamma /\sqrt{n})\times
(0,0,1,1,0,0,0,0)^{\prime }$ for 101 equidistant values of $\gamma $ between 
$0$ and $8$. The worst-case behavior in this setup is qualitatively similar
to the one in Setup~I but slightly less pronounced; we do not report the
results here for brevity. In Setup~III, we again perform the same Monte
Carlo study as in Setup~I, but now with $\theta _{n}=\theta _{0}+(\gamma /%
\sqrt{n})\times (0,0,1,1,0,1/10,1/10,1/10)^{\prime }$ for 101 equidistant
values of $\gamma $ between $0$ and $80$. Note that here we consider a range
for $\gamma $ wider than that in Scenario I and II, where we had $0\leq
\gamma \leq 8$. Figure~2 gives the results for Setup III.

\begin{center}
\begin{tabular}{cc}
\epsfxsize=7cm\epsfbox{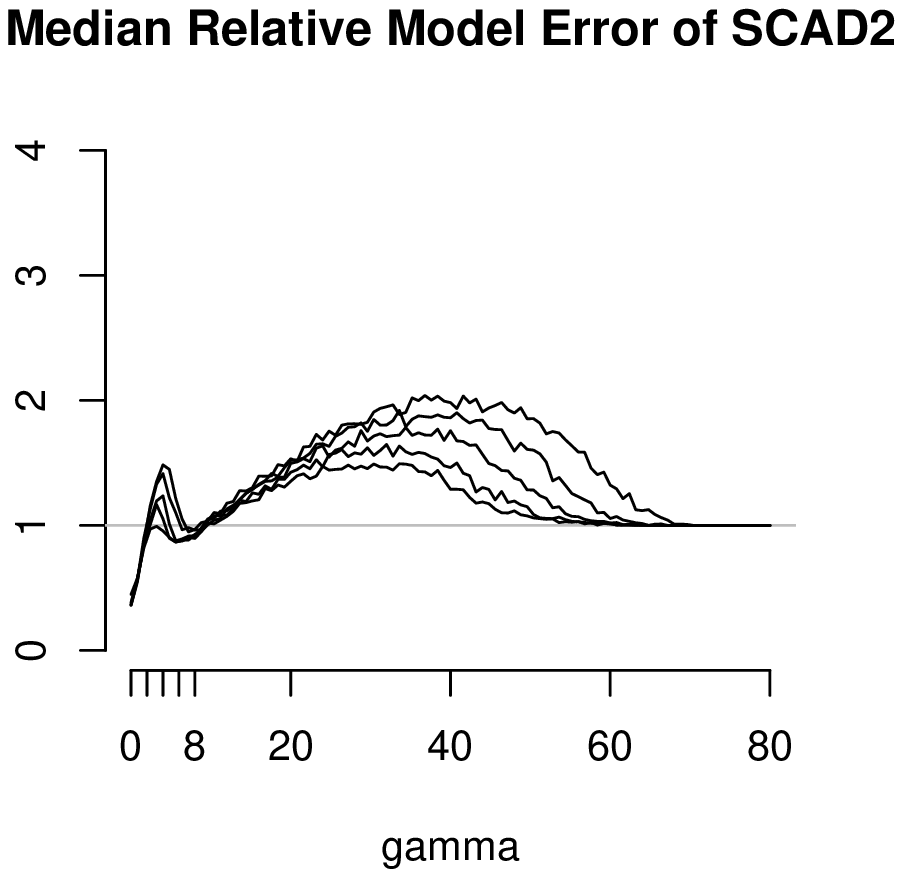} & \epsfxsize=7cm\epsfbox{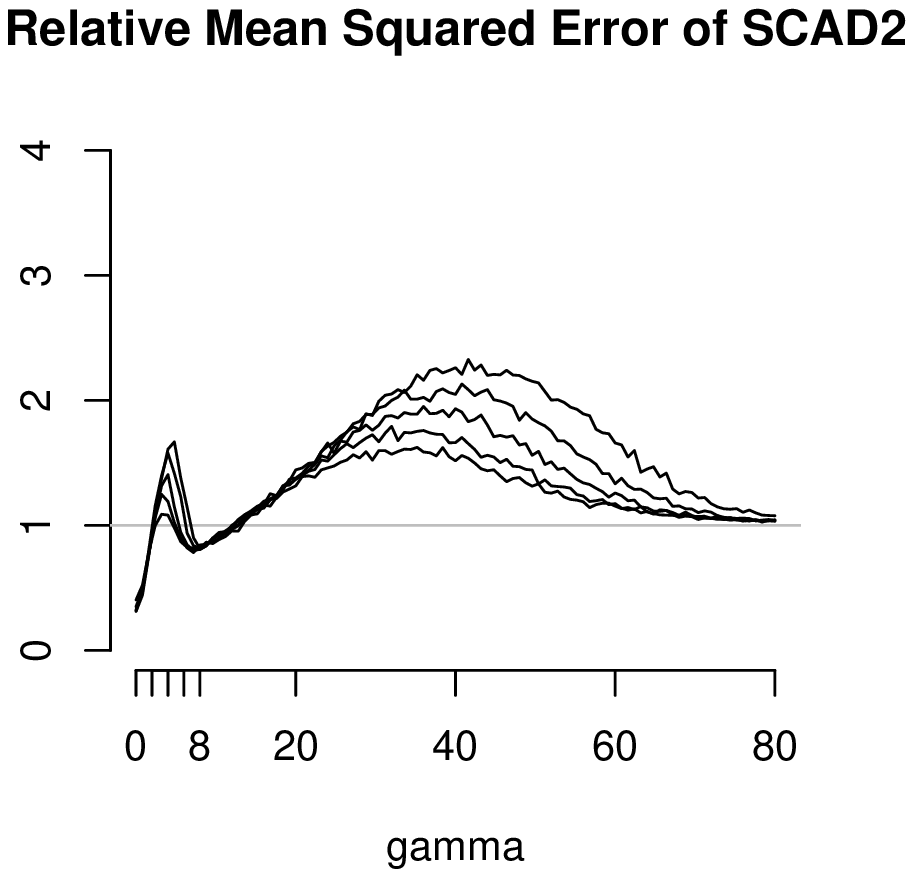}%
\end{tabular}
\end{center}

\begin{quote}
{\small Figure~2: Monte Carlo performance estimates under the true parameter 
$\theta_n = \theta_0 + (\gamma/\sqrt{n})\times
(0,0,1,1,0,1/10,1/10,1/10)^{\prime}$, as a function of $\gamma$. See the
legend of Figure~1 for a description of the graphics. }
\end{quote}

The same considerations as given for Figure~1 also apply to Figure~2. The
new feature in Figure~2 is that the curves are bimodal. Apparently, this is
because now there are two regions, in the range of $\gamma$'s under
consideration, for which some components of the underlying regression
parameter $\theta_n$ are neither very close to zero nor quite large
(relative to sample size): Components $3$ and $4$ for $\gamma$ around $5$
(first peak), and components $6$, $7$, and $8$ for $\gamma$ around $40$
(second peak).

\noindent \textbf{Setups~IV and~V:} Here we perform the same simulations as
in Setup~I, but now with the range of $\lambda $'s considered for
generalized cross-validation given by $\{\delta (\hat{\sigma}/\sqrt{n}%
)(n/60)^{1/10}:\;\delta =0.9,1.1,1.3,\dots ,2\}$ for Setup~IV, and by $%
\{\delta (\hat{\sigma}/\sqrt{n})(n/60)^{1/4}:\;\delta =0.9,1.1,1.3,\dots
,2\} $ for Setup~V. Setup~IV gives `less sparse' estimates while Setup~V
gives `more sparse' estimates relative to Setup~I. The results are
summarized in Figures~3 and~4 below. Choosing the SCAD tuning-parameter $%
\lambda $ so that the estimator is `more sparse' clearly has a detrimental
effect on the estimator's worst-case performance.

\begin{center}
\begin{tabular}{cc}
\epsfxsize=7cm\epsfbox{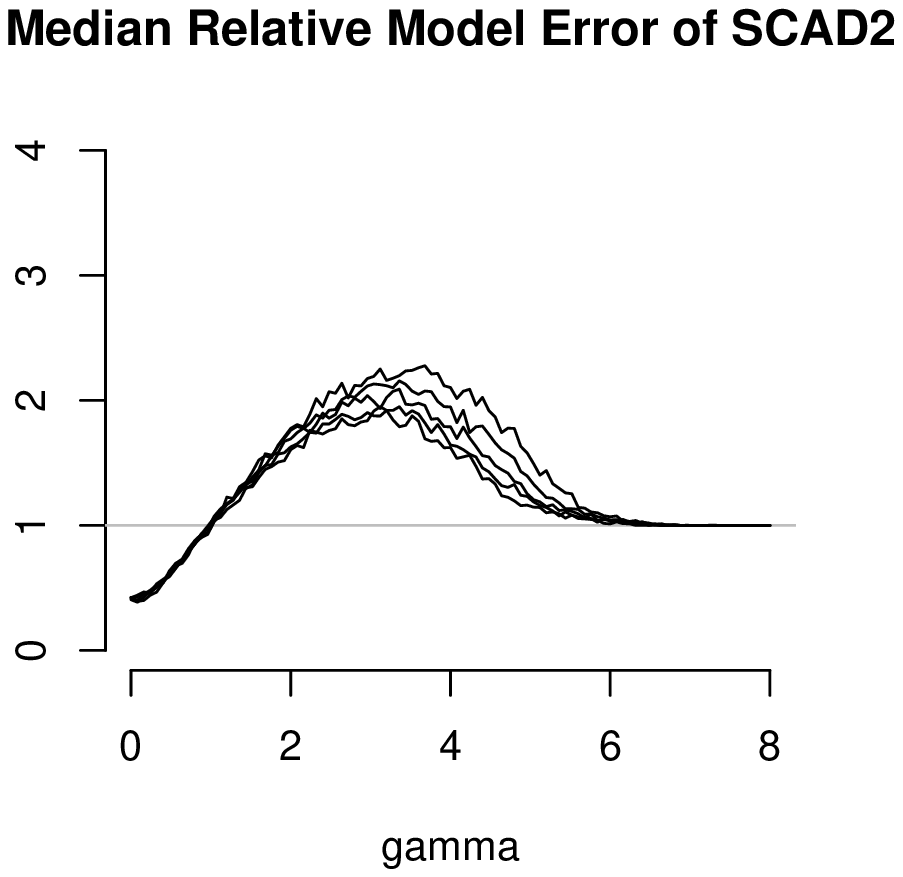} & \epsfxsize=7cm\epsfbox{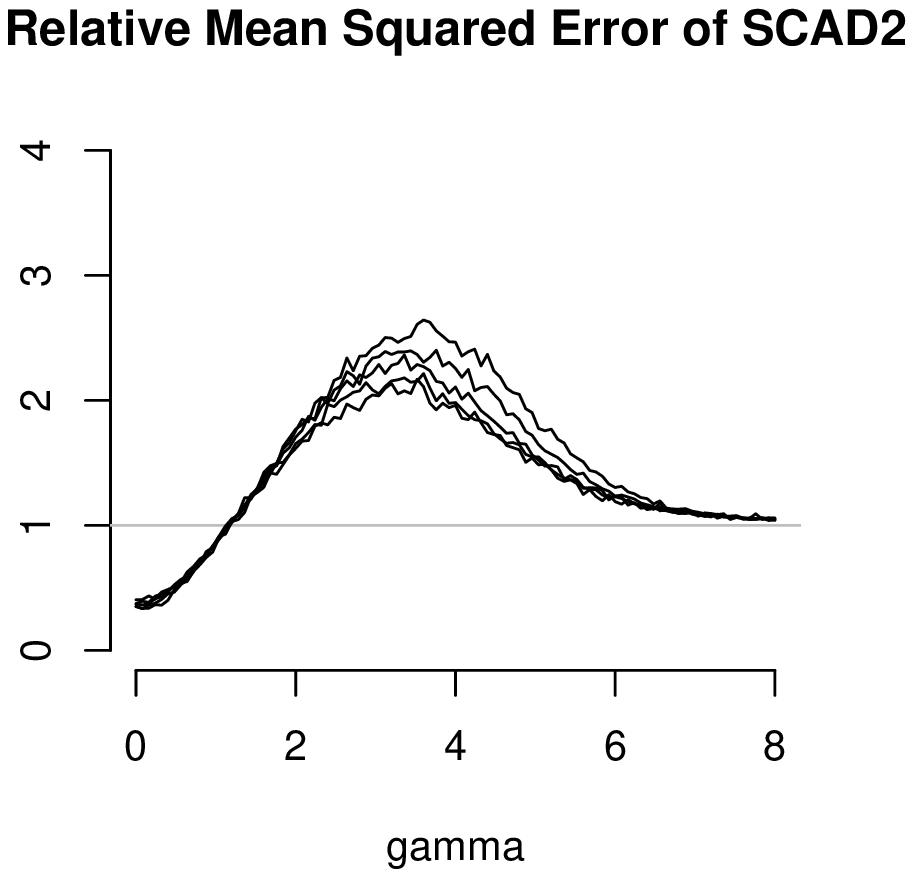}%
\end{tabular}
\end{center}

\begin{quote}
{\small Figure~3: Monte Carlo performance estimates under the true parameter 
$\theta_n = \theta_0 + (\gamma/\sqrt{n})\times (0,0,1,1,0,1,1,1)^{\prime}$
as a function of $\gamma$; the SCAD tuning parameter $\lambda$ is chosen as
described in Setup~IV. }
\end{quote}

\begin{center}
\begin{tabular}{cc}
\epsfxsize=7cm\epsfbox{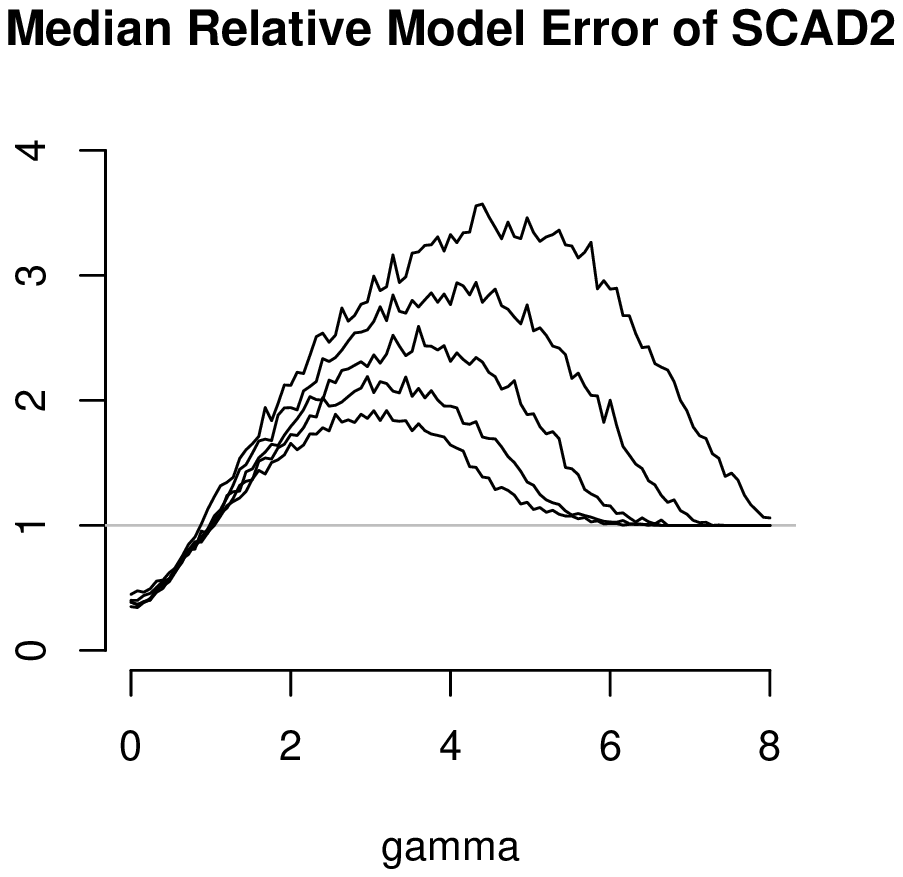} & \epsfxsize=7cm\epsfbox{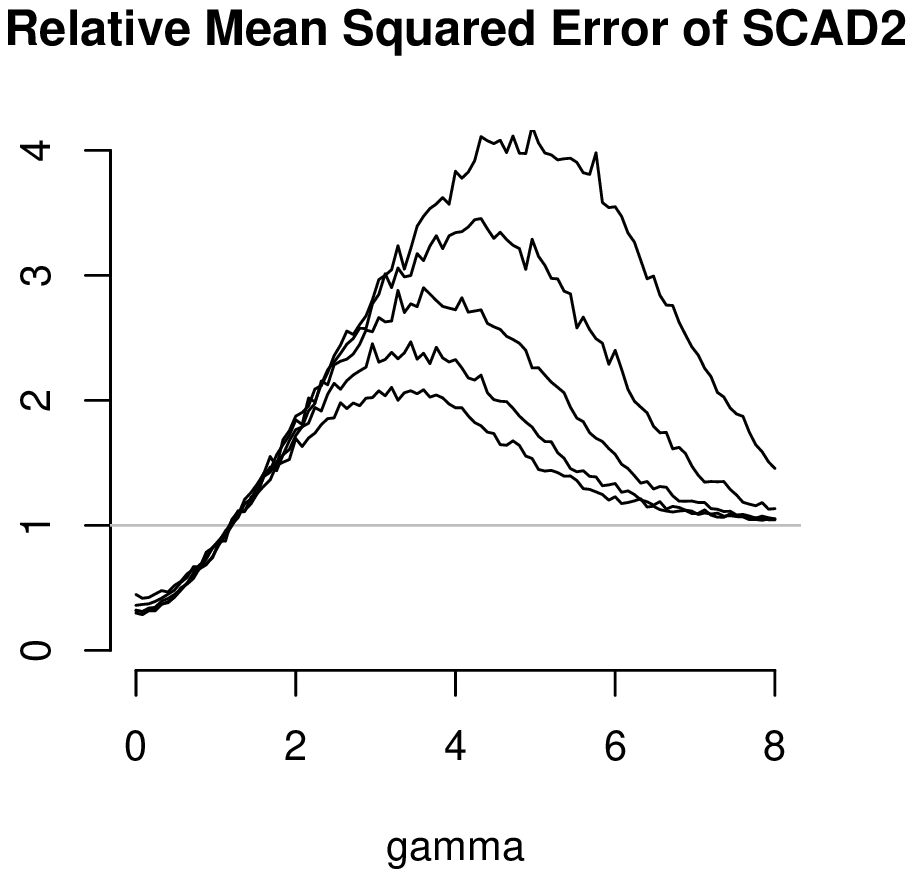}%
\end{tabular}
\end{center}

\begin{quote}
{\small Figure~4: Monte Carlo performance estimates under the true parameter 
$\theta _{n}=\theta _{0}+(\gamma /\sqrt{n})\times (0,0,1,1,0,1,1,1)^{\prime
} $, as a function of $\gamma $; the SCAD tuning parameter $\lambda $ is
chosen as described in Setup~V. }
\end{quote}

In all setups considered so far we have enforced the conditions $\lambda
\rightarrow 0$ and $\sqrt{n}\lambda \rightarrow \infty $ to guarantee
sparsity of the resulting SCAD estimator as risk properties of sparse
estimators are the topic of the paper. In response to a referee we further
consider Setup VI which is identical to Setup I, except that the range of $%
\lambda $'s over which generalized cross-validation is effected is given by $%
\{\delta (\hat{\sigma}/\sqrt{n}):\;\delta =0.9,1.1,1.3,\dots ,2\}$, which is
precisely the range considered in Fan and Li (2001). Note that the resulting 
$\lambda $ does now \emph{not} satisfy the conditions for sparsity given in
Theorem 2 of Fan and Li (2001). The results are shown in Figure 5 below. The
findings are similar to the results from Setup I, in that SCAD2 gains over
the least squares estimator in a neighborhood of $\theta _{0}$, but is worse
by approximately a factor of $2$ over considerable portions of the range of $%
\gamma $, showing once more that the simulation study in Fan and Li (2001)
does not tell the entire truth. What is, however, different here from the
results obtained under Setup I is that -- not surprisingly at all -- the
worst case behavior now does not get worse with increasing sample size.
[This is akin to the boundedness of the worst case risk of a
post-model-selection estimator based on a conservative model selection
procedure like AIC or pre-testing with a sample-size independent critical
value.]

\begin{center}
\begin{tabular}{cc}
\epsfxsize=7cm\epsfbox{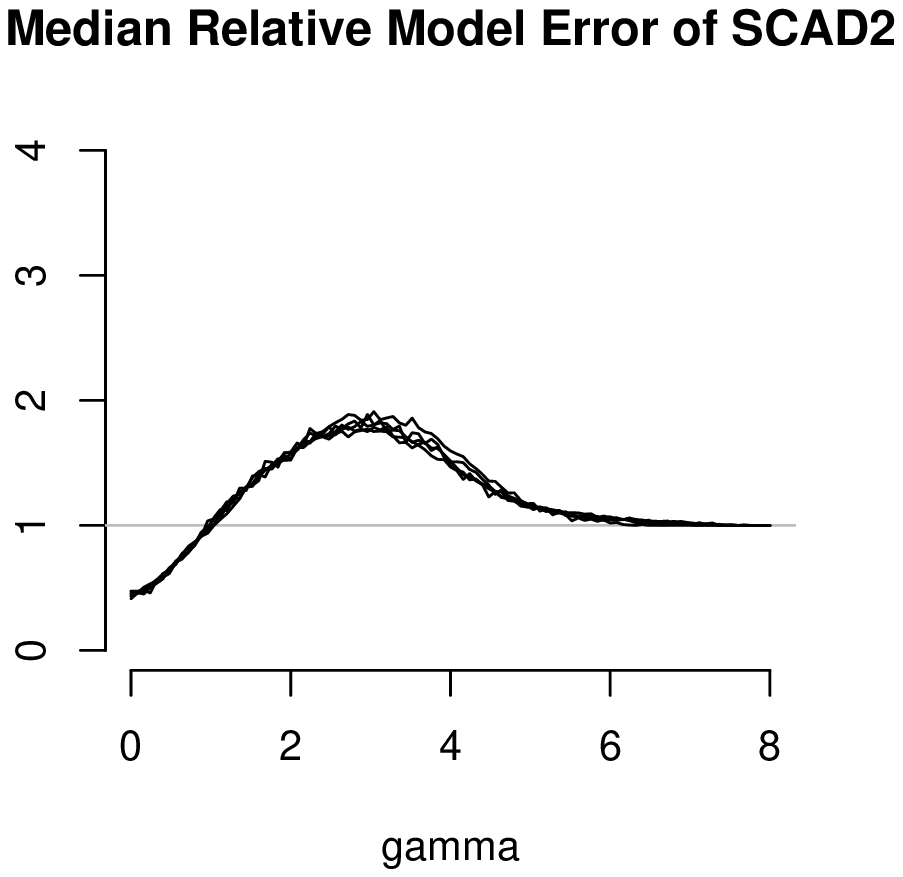} & \epsfxsize=7cm%
\epsfbox{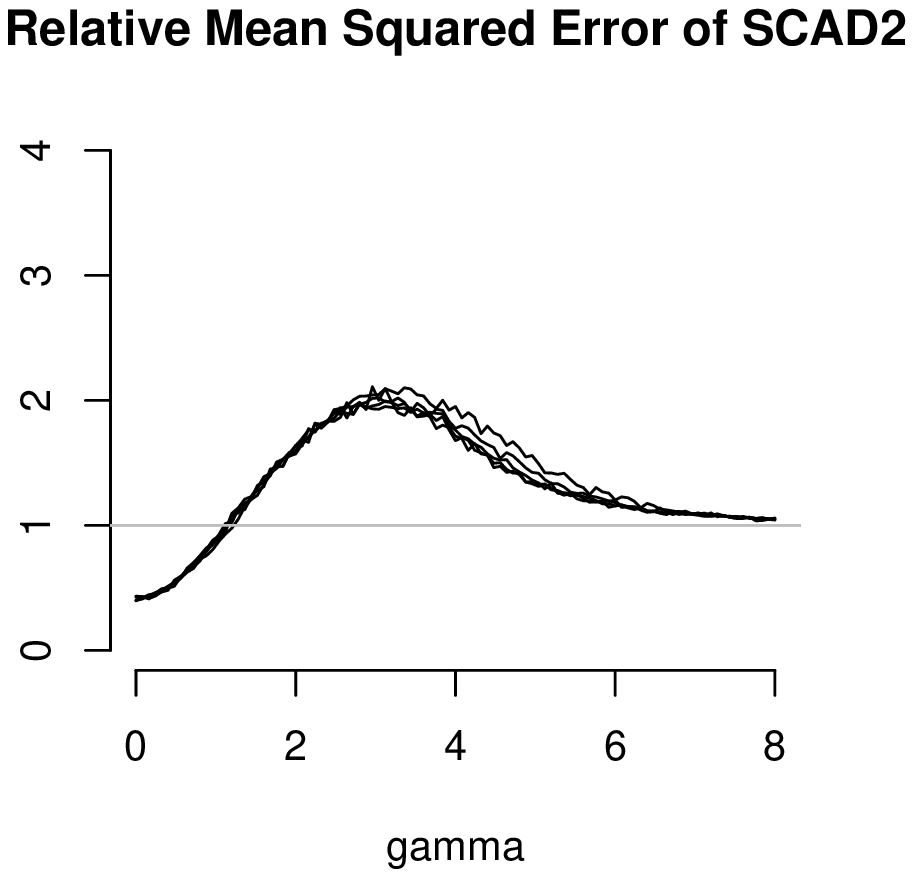}%
\end{tabular}
\end{center}

\begin{quote}
{\small Figure~5: Monte Carlo performance estimates under the true parameter 
$\theta _{n}=\theta _{0}+(\gamma /\sqrt{n})\times (0,0,1,1,0,1,1,1)^{\prime
} $, as a function of $\gamma $; the SCAD tuning parameter $\lambda $ is
chosen as described in Setup~VI.}
\end{quote}

\section{Conclusion}

We have shown that sparsity of an estimator leads to undesirable risk
properties of that estimator. The result is set in a linear model framework,
but easily extends to much more general parametric and semiparametric
models, including time series models. Sparsity is often connected to a
so-called \textquotedblleft oracle property\textquotedblright . We point out
that this latter property is highly misleading and should not be relied on
when judging performance of an estimator. Both observations are not really
new, but worth recalling: Hodges' construction of an estimator exhibiting a
deceiving pointwise asymptotic behavior (i.e., the oracle property in
today's parlance) has led mathematical statisticians to realize the
importance uniformity has to play in asymptotic statistical results. It is
thus remarkable that today -- more than 50 years later -- we observe a
return of Hodges' estimator in the guise of newly proposed estimators (i.e.,
sparse estimators). What is even more surprising is that the deceiving
pointwise asymptotic properties of these estimators (i.e., the oracle
property) are now advertised as virtues of these methods. It is therefore
perhaps fitting to repeat Hajek's (1971, p.153) warning:%
\begin{eqnarray*}
&&\text{\textquotedblleft Especially misinformative can be those limit
results that are not uniform. Then the limit } \\
&&\text{may exhibit some features that are not even approximately true for
any finite }n\text{.\textquotedblright }
\end{eqnarray*}%
The discussion in the present paper as well as in Leeb and P\"{o}tscher
(2005) shows in particular that distributional or risk behavior of
consistent post-model-selection estimators is not as sometimes believed, but
is much worse.

The results of this paper should not be construed as a criticism of
shrinkage-type estimators including penalized least squares (maximum
likelihood) estimators per se. Especially if the dimension of the model is
large relative to sample size, some sort of shrinkage will typically be
beneficial. However, achieving this shrinkage through sparsity is perhaps
not such a good idea (at least when estimator risk is of concern). It
certainly cannot simply be justified through an appeal to the oracle
property.\footnote{%
In this context we note that \textquotedblleft
superefficiency\textquotedblright\ per se is not necessarily detrimental in
higher dimensions as witnessed by the Stein phenomenon. However, not all
forms of \textquotedblleft superefficiency\textquotedblright\ are created
equal, and \textquotedblleft superefficiency\textquotedblright\ generated
through sparsity of an estimator typically belongs to the undesirable
variety as shown in the paper.}

\section*{Acknowledgements}

A version of this paper was previously circulated in 2004. We are grateful
to the editor Ron Gallant and the referees as well as to Hemant Ishwaran,
Paul Kabaila, Richard Nickl, and Yuhong Yang for helpful comments.

\bibliographystyle{econometrica}
\bibliography{lit}

\section{References}

\quad \thinspace \thinspace Bunea, F. (2004): Consistent covariate selection
and post model selection inference in semiparametric regression. \emph{%
Annals of Statistics }32, 898-927.

Bunea, F. \& I.~W. McKeague (2005): Covariate selection for semiparametric
hazard function regression models. \emph{Journal of Multivariate Analysis }%
92, 186-204.

Cai, J., Fan, J., Li, R., \& H. Zhou (2005): Variable selection for
multivariate failure time data, \emph{Biometrika} 92, 303-316.

Fan, J. \& R. Li (2001): Variable selection via nonconcave penalized
likelihood and its oracle properties. \emph{Journal of the American
Statistical Association} 96, 1348-1360.

Fan, J. \& R. Li (2002): Variable selection for Cox's proportional hazards
model and frailty model. \emph{Annals of Statistics }30, 74-99.

Fan, J. \& R. Li (2004): New estimation and model selection procedures for
semiparametric modeling in longitudinal data analysis. \emph{Journal of the
American Statistical Association} 99, 710-723.

Fan, J. \& H. Peng (2004): Nonconcave penalized likelihood with a diverging
number of parameters. \emph{Annals of Statistics }32, 928-961.

Foster D.~P. \& E.~I. George (1994): The risk inflation criterion for
multiple regression. \emph{Annals of Statistics }22, 1947-1975.

Frank, I.~E. \& J.~H. Friedman (1993): A statistical view of some
chemometrics regression tools (with discussion). \emph{Technometrics \ }35,
109-148.

Hajek, J. (1971): Limiting properties of likelihoods and inference. In:
V.~P. Godambe and D.~A. Sprott (eds.), \emph{Foundations of Statistical
Inference: Proceedings of the Symposium on the Foundations of Statistical} $%
\emph{Inference}$\emph{, University of Waterloo, Ontario, March 31 -- April
9, 1970}, 142-159. Toronto: Holt, Rinehart \& Winston.

Hajek, J. \& Z. Sidak (1967): \emph{Theory of Rank Tests}. New York:
Academic Press.

Hosoya, Y. (1984): Information criteria and tests for time series models.
In: O.~D. Anderson (ed.), \emph{Time Series Analysis: Theory and Practice \ }%
5, 39-52. Amsterdam: North-Holland.

Judge, G.~G. \& M.~E. Bock (1978): \emph{The Statistical Implications of
Pre-test and Stein-rule Estimators in Econometrics}. Amsterdam:
North-Holland.

Kabaila, P. (1995): The effect of model selection on confidence regions and
prediction regions. \emph{Econometric Theory} {\ 11}, 537-549.

Kabaila, P. (2002): On variable selection in linear regression. \emph{%
Econometric Theory} {\ 18}, 913-915.

Knight, K. \& W. Fu (2000): Asymptotics of lasso-type estimators. \emph{%
Annals of Statistics \ }28, 1356-1378.

Koul, H.~L. \& W. Wang (1984): Local asymptotic normality of randomly
censored linear regression model. \emph{Statistics \& Decisions},\emph{\ }%
Supplement Issue\emph{\ }No.\emph{\ }1, 17-30.

Lehmann, E.~L. \& G. Casella (1998): \emph{Theory of Point Estimation}.
Springer Texts in Statistics. New York: Springer-Verlag.

Leeb, H. \& B.~M. P\"{o}tscher (2005): Model selection and inference: facts
and fiction. \emph{Econometric Theory} {\ 21}, 21-59.

Leeb, H. \& B.~M. P\"{o}tscher (2006): Performance limits for estimators of
the risk or distribution of shrinkage-type estimators, and some general
lower risk-bound results. \emph{Econometric Theory} {\ 22}, 69-97.
(Correction, ibid., forthcoming.)

P\"{o}tscher, B.~M. (1991): Effects of model selection on inference. \emph{%
Econometric Theory} {\ 7}, 163-185.

Shibata R. (1986a): Consistency of model selection and parameter estimation. 
\emph{Journal of Applied Probability, Special Volume \ }23A, 127-141.

Shibata R. (1986b): Selection of the number of regression variables; a
minimax choice of generalized FPE. \emph{Annals of the Institute of
Statistical Mathematics \ }38, 459-474.

Tibshirani, R.~J. (1996): Regression shrinkage and selection via the LASSO. 
\emph{Journal of the Royal Statistical Society, Ser. B \ }58, 267-288.

Von Rosen, D. (1988): Moments for the inverted Wishart distribution. \emph{%
Scandinavian Journal of Statistics \ }15, 97-109.

Yang, Y. (2005): Can the strengths of AIC and BIC be shared? A conflict
between model identification and regression estimation. \emph{Biometrika \ }%
92, 937-950.

Zou, H. (2006): The adaptive lasso and its orcale properties. \emph{Journal
of the American Statistical Association }101, 1418-1429.

\end{document}